\documentclass[12pt]{article}
\usepackage{amssymb,latexsym,amsmath}
\def\C{\mathbb C}

\def\N{\mathbb N}

\def\D{\mathbb D}

\newtheorem{thm}{Theorem}[section]
\newtheorem{lem}{Lemma}[section]

\begin{document}
\sffamily

\title { Wiman-Valiron theory for a class of functions meromorphic in the unit disc}
\author{J.K. Langley and John Rossi}
\maketitle

\begin{abstract}
Analogues of the key results of  Wiman-Valiron theory are proved for a class of functions meromorphic in the unit disc, 
based on an approach developed by Bergweiler, Rippon and Stallard for the plane setting. The results give local 
approximations for the function and its logarithmic derivative and, in the case of positive order of growth, for higher
order logarithmic derivatives as well. 
\end{abstract}

\section{Introduction}

\noindent
Classical Wiman-Valiron theory describes the behaviour of an entire function $f(z)$  
by analyzing its power series $\sum_{n=0}^{\infty} a_n z^n$ (see \cite{Hayman} for the definitive reference). 
For $r>0$ we define the {\em maximum term} $\mu(r) = \max_{n\geq 0}|a_n|r^n$, and  the {\em central index} $N(r)$ 
is then the largest integer $n$ for which this maximum is attained. 
A seminal result  of the theory states that if $ \gamma > 1/2$ and $M \in \N$ and
$r \in [1, \infty) \setminus E$, where $E\subseteq [1, \infty)$ is a set of finite logarithmic measure, and if $|f(z_r)| =
M(r, f) =\max_{|\zeta|=r} |f(\zeta)|$, then
\begin{equation}\label{asym} f(z) \sim \left(\frac{z}{z_r}\right)^{N(r)} f(z_r) \end{equation}
and 
\begin{equation}\label{asym1}f^{(q)}(z) \sim \left(\frac{N(r)}{z}\right)^qf(z)\end{equation}
for $| \log (z/z_r)| < N(r)^{-\gamma}$ and $1 \leq q \leq M$. 
Equations (\ref{asym}) and (\ref{asym1}) imply that near maximum modulus points, 
$f$ behaves like  a monomial, namely the dominant term of its power series. 
This has proved decisive in numerous applications, including to differential equations \cite{Lai1} and iteration theory \cite{Ber1,Ber3,Er0}.

Two recent papers have generalized these results in different ways. Fenton and Rossi \cite[Theorem 1]{FR} used the power series
approach to obtain the 
approximation (\ref{asym1}) at points where $f$ is close to its maximum modulus,
when $f$ is a function analytic in the unit disc of positive order
of growth (as defined by (\ref{a1}) below with $B(r,v) = \log M(r,f)$).
Bergweiler, Rippon and Stallard \cite[Theorem 2.2]{BRS} developed a powerful technique, not involving power series and
closer in spirit to work of Macintyre \cite{Mac}, resulting in a Wiman-Valiron theory for certain classes of 
plane meromorphic functions. Their results in \cite{BRS} (and the extensions by Bergweiler
in \cite{bergsizediscs}) are applicable, in particular, to any transcendental
meromorphic function in the plane for which the inverse function has a direct singularity over infinity \cite{BRS}.  
With these results in mind, it seems  natural to seek a result analogous to 
that of \cite{BRS} for the unit disc.  Before stating our theorems we need some definitions.

Let the function $f$ be meromorphic in the unit disc $\D = D(0, 1) = \{ z \in \C : |z| < 1 \}$. In analogy with \cite{BRS},
a \textit{direct tract} of $f$ will mean a component $U$ 
of the set $\{ z \in \D : \,  |f(z)| > R \}$, for some $R \in (0, \infty)$, such that $U$ contains no poles of $f$ but $f$ is unbounded 
on $U$. It follows using the maximum principle that $\partial U$ meets the circle  $|z| = 1$. 

Assume henceforth that $f$ has a direct tract: then as in \cite{BRS} the function 
\begin{equation}
 \label{1}
v(z) = \log \left| \frac{f(z)}{R} \right| \quad (z \in U), \quad v(z) = 0 \quad (z \in \D \setminus U) ,
\end{equation}
is continuous, subharmonic and unbounded on $\D$. 
For $0 < r < 1$ let 
\begin{equation}
 \label{2}
B(r) = B(r, v) = \max \{ v(z) : |z| = r \} , \quad a(r) = r B'(r) = \frac{dB(r)}{d \log r } .
\end{equation}
Here $B(r)$ is a non-decreasing convex function of $\log r$ for $0 < r < 1$,
and $a(r)$ (which is taken to be the right derivative with respect to $\log r$  at those countably many 
points at which $B$ is not differentiable) is non-decreasing and tends to $\infty$ as $r \to 1-$. 
Choose $r_0 \in (0, 1)$ and $ \beta, \delta $ with 
\begin{equation}
 \label{3}
\hbox{$B(r) \geq 2$ and $a(r) \geq 2 $ for $r_0 \leq r < 1$},  \quad  0  < \beta \leq \frac12 , \quad  \delta > 0.
\end{equation}

\begin{thm}
 \label{thm1}
Let the function $f$ be meromorphic in the unit disc $\D$ with a direct tract. Using  the notation 
(\ref{1}), (\ref{2}) and (\ref{3}) set 
\begin{equation}
 \label{3a}
\varepsilon (r) = \min \left\{ \frac{1-r}{2 a(r)^\beta (\log a(r))^{1 +  \delta } }, \quad  
\frac{ 1}{  a(r)^{1-\beta}(\log a(r))^{1 +  \delta } }\right\} 
\end{equation}
for $r_0 \leq r < 1$. Then there exists a set $E \subseteq [r_0, 1)$ satisfying
\begin{equation}
 \label{4}
\int_E \, \frac{dt}{1-t} < \infty ,
\end{equation}
such that, as $r \to 1-$ with $r \not \in  E$, if $z_r$ is chosen with
$|z_r| = r$ and $v(z_r) = B(r, v)$ then 
\begin{equation}
 \label{con1}
f(z) \sim f(z_r) \left( \frac{z}{z_r} \right)^{a(r)} = f(z_r) \exp \left( a(r) \log \frac{z}{z_r} \right) 
\end{equation}
and 
 \begin{equation}
 \label{20}
\frac{f'(z)}{f(z)} \sim \frac{ a(r)}z  
\end{equation}
for $| z- z_r | < \varepsilon (r)/2048$.
\end{thm}
Here the logarithm in (\ref{con1}) is chosen so as to vanish at $z_r$. In particular, Theorem \ref{thm1} certainly applies if $f$ is analytic and 
unbounded in the unit disc, and gives analogues of (\ref{asym}) for $f$ and (\ref{asym1}) for $q=1$. However, 
it cannot be expected that an  estimate 
\begin{equation}
 \label{b1}
L_q(z) = \frac{f^{(q)}(z)}{f(z)} \sim \frac{a(r)^q}{z^q} 
\end{equation}
always holds for $q \geq 2$, even at $z_r$ itself, as is shown by the  well known example
\begin{equation}
 \label{exa1}
f(z) = (1-z)^{-\gamma} , \quad \gamma > 0, 
\end{equation}
for which, with $R=1$,  
$$B(r) = \gamma \log \frac1{1-r} \quad \hbox{and} \quad a(r) = \frac{\gamma r }{1-r},$$
but $f''(z)/f(z)$ is never asymptotic to $(f'(z)/f(z))^2$. On the other hand, 
if the function $v$ in (\ref{1}) has positive order
\begin{equation}
 \label{a1}
\rho(v) = \limsup_{r \to 1-} \frac{\log B(r, v)}{\log 1/(1-r) } > \rho_0 > 0 ,
\end{equation}
then more can be proved. It follows easily from (\ref{2}) and (\ref{a1}) that 
\begin{equation*}
\limsup_{r \to 1-} \frac{\log a(r)}{\log 1/(1-r) } > 1+ \rho_0 > 1, 
\end{equation*}
and if $\beta$ is chosen small enough in (\ref{3}) then 
\begin{equation}
 \label{a1a}
\limsup_{r \to 1-} \frac{\log a(r)}{\log 1/(1-r) } > \frac{1+ \rho_0 }{1 - 2 \beta} .  
\end{equation}
Hence there exists a sequence $(r_n)$ satisfying
\begin{equation}
 \label{a1b}
r_0 \leq r_n < 1, \quad \lim_{n \to \infty} r_n = 1, \quad 
\lim_{n \to \infty} (1-r_n)^{1+\rho_0} a(r_n)^{1-2 \beta} = \infty .
\end{equation}
The following theorem will be proved. 

\begin{thm}
 \label{thm2}
Let the function $f$ be meromorphic in the unit disc $\D$ with a direct tract and, using  the notation 
(\ref{1}), (\ref{2}) and (\ref{3}),
assume that $v$ and $\beta$ satisfy (\ref{a1}) and (\ref{a1a}). Let the set $E$ be as in Theorem~\ref{thm1},
let  $(r_n)$ be any sequence satisfying (\ref{a1b}), and let $M$ be a positive integer. Then for all sufficiently large $n$ 
and all $r$ satisfying
\begin{equation}
 \label{a1c}
r_n \leq r \leq r_n' = 1- (1-r_n)^{1 + \rho_0} , \quad r \not \in E, 
\end{equation}
the function $f$ satisfies (\ref{con1}) and (\ref{b1}) for $1 \leq q \leq M$ and 
\begin{equation}
 \label{a1d}
|z-z_r | < \frac{\varepsilon(r)}{2048} = \frac{ 1}{ 2048  a(r)^{1-\beta}(\log a(r))^{1 +  \delta } } ,
\end{equation}
where $z_r$ is as in Theorem \ref{thm1}. 
\end{thm}
Here 
$$
\int_{r_n}^{r_n'} \, \frac{dt}{1-t} = \rho_0 \log \frac1{1-r_n}
$$
and so in view of (\ref{4}) the set of $r$ satisfying (\ref{a1c}) comprises most of the interval $[r_n, r_n']$. 
Moreover, if $\rho(v)$ is large enough then $\beta$ may be chosen close to $1/2$ in
(\ref{3}), (\ref{3a}) and (\ref{a1a}). 

\section{A growth lemma}

\begin{lem}
 \label{growthlem}
Let $x_0$ and $\delta$ be positive,  let $0 < \beta \leq 1/2$, and let 
$A: [x_0, \infty) \to (1, \infty)$  be a non-decreasing  function. Then there exists a set $E_0 \subseteq [x_0, \infty)$, of 
finite linear measure, such that, for $x \in [x_0, \infty) \setminus E_0$, 
\begin{eqnarray}
 \label{10}
A \left( x + \frac1{A(x)^\beta (\log A(x))^{1 + \delta}} \right)  &<&  A(x) + A(x)^{1-\beta} \nonumber \\
\quad \hbox{and}  \quad A \left( x - \frac1{A(x)^\beta (\log A(x))^{1 + \delta}} \right) &>& A(x) - A(x)^{1-\beta} . 
\end{eqnarray}
\end{lem}
\textit{Proof.} This follows at once from \cite[Lemma 2.1]{bergsizediscs} with the choice
$\sigma_2(t) = t^{1-\beta} $, $\sigma_1(t) = t^\beta (\log t)^{1+\delta} $. 
\hfill$\Box$
\vspace{.1in}

\section{Proof of Theorem \ref{thm1}}

Let $f$ be as in the hypotheses of Theorem \ref{thm1}, and
denote by $C$ positive constants, not necessarily the same at each occurrence, but always independent of $r$. 
The set $E$ is determined by the following lemma. 

\begin{lem}
 \label{lemexset}
There exists a set $E \subseteq [r_0, 1)$ satisfying (\ref{4}) 
such that, for $r \in [r_0, 1) \setminus E$, 
\begin{equation}
 \label{5}
a(r + \varepsilon(r)) < a(r) + a(r)^{1-\beta} 
\end{equation}
and \begin{equation}
 \label{6}
a(r - \varepsilon(r)) > a(r) - a(r)^{1-\beta} ,
\end{equation}
as well as 
\begin{equation}
 \label{7}
(1-r) a(r) < B(r)^{1 + \beta } .
\end{equation}
\end{lem}
\textit{Proof.} To establish (\ref{5}) and (\ref{6}) set
\begin{equation}
 \label{8}
x = x(r) = \log \frac1{1-r}, \quad A(x) = a(r), \quad x_0 = \log \frac1{1-r_0} \leq x < \infty .
\end{equation}
By Lemma \ref{growthlem} there exists a set $F_1 \subseteq [x_0, \infty) $ with
\begin{equation}
 \label{9}
\infty > \int_{F_1} \, dx = \int_{E_1} \, \frac{dr}{1-r} , \quad 
E_1 = \{ r \in [r_0, 1) : x(r) \in F_1 \} ,
\end{equation}
such that (\ref{10}) holds for $x \in [x_0, \infty) \setminus F_1$.  
For $r \in [r_0, 1) \setminus E_1$ and $x = x(r) $ define $r'$ and $r''$ by
$$
\log \frac1{1-r'} = x + D(r) , \quad
\log \frac1{1-r''} = x -  D(r), 
$$
in which
$$  D(r) = \frac1{A(x)^\beta (\log A(x))^{1 + \delta} } = \frac1{a(r)^\beta (\log a(r))^{1 + \delta} } .
$$
Then 
$$
1-r' = (1-r) e^{-D(r)} , \quad 1-r'' = (1-r) e^{D(r)}
$$
and so, as $r \to 1-$, by (\ref{3a}),  
$$
r'-r = (1-r)D(r) (1+ o(1)) \geq \varepsilon (r) , \quad
r-r'' = (1-r)D(r) (1+ o(1)) \geq \varepsilon (r) ,
$$
which gives (\ref{5}) and (\ref{6}), using (\ref{10}) and the fact that $a(r)$ is non-decreasing. 

Next, let 
\begin{equation}
 \label{12}
E_2 = \{ r \in [r_0, 1) : (1-r) a(r) \geq B(r)^{1 + \beta } \} .
\end{equation}
Then 
\begin{equation}
 \label{13}
\int_{E_2} \, \frac{dr}{1-r} \leq 
\int_{E_2} \, \frac{a(r)}{B(r)^{1 + \beta }} \, \frac{dr}{r} \leq 
\int_{[r_0, 1)} \, \frac{a(r)}{B(r)^{1 + \beta }} \, \frac{dr}{r} =
\int_{[r_0, 1)} \, \frac{B'(r)}{B(r)^{1 + \beta }} \, dr 
< \infty .
\end{equation}
The proof of the lemma is completed by taking $E = [r_0, r_0'] \cup E_1 \cup E_2$ for some $r_0' \in (r_0, 1)$, and
(\ref{4}) follows from (\ref{9}), (\ref{12}) and (\ref{13}).
\hfill$\Box$
\vspace{.1in} 

\begin{lem}
 \label{lemB(r)}
For $r \in [r_0, 1) \setminus E$ the function $B(r)$ satisfies 
\begin{equation}
 \label{14}
B(s) \leq B(r) + a(r) \log \frac{s}r + \phi(r) \quad \hbox{for} \quad 
r - \varepsilon (r) \leq s \leq r + \varepsilon (r) ,
\end{equation}
in which 
\begin{equation}
 \label{15}
0 \leq \phi(r) \leq C a(r)^{1-\beta}  \varepsilon(r) = o(1) 
\end{equation}
as $r \to 1-$.  
\end{lem}
\textit{Proof.} 
Let $r \in [r_0, 1) \setminus E$ be close to $1$. First take $r \leq s \leq r + \varepsilon (r)$; then (\ref{5}) yields
\begin{eqnarray*}
 B(s) &=& B(r) + \int_r^s \, a(t) \, \frac{dt}{t} \\
&\leq& B(r) + \int_r^s (a(r) + a(r)^{1-\beta} ) \, \frac{dt}{t} \\
&\leq& B(r) + a(r) \log \frac{s}{r} + a(r)^{1-\beta} \log (1 + \varepsilon(r)/r) \\
&\leq& B(r) + a(r) \log \frac{s}{r} + C a(r)^{1-\beta}  \varepsilon(r). 
\end{eqnarray*}
Similarly, $r - \varepsilon (r) \leq s \leq r$ and (\ref{6}) give
\begin{eqnarray*}
 B(s) &=& B(r) - \int_s^r \, a(t) \, \frac{dt}{t} \\
&\leq& B(r) - \int_s^r  (a(r) - a(r)^{1-\beta} ) \, \frac{dt}{t} \\
&\leq& B(r) + a(r) \log \frac{s}{r} +   a(r)^{1-\beta}  \log \frac{r}{s}  \\
&\leq& B(r) + a(r) \log \frac{s}{r} + a(r)^{1-\beta} \log \frac1{1 - \varepsilon(r)/r} \\
&\leq& B(r) + a(r) \log \frac{s}{r} + C a(r)^{1-\beta}  \varepsilon(r). 
\end{eqnarray*}
In view of  (\ref{3a}), the lemma follows. 
\hfill$\Box$
\vspace{.1in}

\begin{lem}
 \label{disclem}
Let $r \in [r_0, 1) \setminus E$, set $\sigma = \sigma(r) = \varepsilon (r)/2048$ and choose $z_r $
with $|z_r| = r$ and $v(z_r) = B(r, v)$. If $r $ is close enough to $1$ then the disc $D(z_r, 4 \sigma) $ of centre $z_r$ and radius
$4 \sigma$ lies in $ U$. 
\end{lem}
\textit{Proof.} In the  argument below the underlying ideas are the same as for
the corresponding lemma in \cite{BRS}, but the method is simplified somewhat insofar as
the Riesz decomposition of a subharmonic function is not required.
Choose $R' $ with $R' - R$ small and positive,
such that $f$ has no critical points $z$ with $|f(z)| = R'$.
Following \cite{BRS} and using (\ref{14}) and (\ref{15}) form the subharmonic function
\begin{equation}
 \label{16}
u(z) = v(z) - B(r) - a(r) \log \frac{|z|}r \leq \phi(r) = o(1) 
\end{equation}
on $D(z_r, 2048 \sigma)$. For $z $ in $ D(z_r, 2048 \sigma)$,
formulas  (\ref{3a}) and (\ref{7}) give
\begin{equation}
 \label{17}
\left| a(r) \log \frac{|z|}r  \right| \leq C a(r) \varepsilon (r) 
\leq \frac{C (1-r) a(r)}{a(r)^\beta (\log a(r))^{1+\delta} }  
\leq \frac{C B(r)^{1+\beta}}{a(r)^\beta (\log a(r))^{1+\delta} }  
= o( B(r)) ,
\end{equation}
using the  inequality $B(r) \leq C a(r) + C$, which follows from integration of $a(t)$ with respect to $\log t$.

Assume that $r \in [r_0, 1) \setminus E$ is close to $1$ and that
the assertion of the lemma is false. 
Let $U'$ be the component of the set
$\{ z \in \D : \, |f(z)| > R' \}$ which contains $z_r$; then there is a component $K$ of $\D \setminus U'$ which meets 
$D(z_r, 4 \sigma) $. Let 
$V$ be the component of $D(z_r, 2048 \sigma) \cap U'$ which contains $z_r$, and let $T$ be the set of 
$t \in (4 \sigma, 1024 \sigma )$ for which the circle 
$|z-z_r| = t$ is contained in $U'$. 

Suppose first that $T$ is empty, 
and set $W = \{ z \in \partial V : \, |z-z_r | = 2048 \sigma \}$.
Then
the standard Carleman-Tsuji estimate for 
harmonic measure \cite[p.112]{Tsuji}  gives
$$
\omega (z_r, W, V) \leq 3 \sqrt{2} \exp \left( - \pi \int_{[4 \sigma, 1024 \sigma ]  } \, \frac{dt}{ 2 \pi t }  \right) 
< \frac12.
$$
Hence the harmonic measure of $\partial V \cap D(z_r, 2048 \sigma)$
with respect to $V$, evaluated at $z_r$, is at least $1/2$, and $u(z) \leq (-1+o(1))B(r)$ for 
$z \in \partial V \cap D(z_r, 2048 \sigma)$, by (\ref{1}),
(\ref{16}) and (\ref{17}). 
Since $u(z_r) = 0$ but
$u(z) \leq o(1)$ on $V$,  by (\ref{16}), 
applying the two-constants theorem gives a contradiction.

It must therefore be the case that $T$ is non-empty. For $0 < t < 2048 \sigma$ set 
$$
I(t) = \frac1{2 \pi} \int_0^{2 \pi} \, u( z_r + t e^{i\theta} ) \, d \theta 
= \frac1{2 \pi} \int_0^{2 \pi} \, v( z_r + t e^{i\theta} ) \, d \theta  - B(r), \quad J(t) = t I'(t) = \frac{dI(t)}{d \log t}  ,
$$
using (\ref{16}) and the mean value property of harmonic functions. 
Here $I(t) \geq u(z_r) = 0$ is a non-decreasing convex function of $\tau = \log t$, 
while $J(t)$ exists for all but countably many $t$ in $(0 ,  2048 \sigma )$ (at these exceptional points one may
take the right derivative), and is also non-decreasing. It will be shown that
\begin{equation}
 \label{atleast1}
J(t) \geq 1 \quad \hbox{ for } \quad 1024 \sigma \leq t < 2048 \sigma . 
\end{equation}
To prove this, let $s \in T$. Then the circle $|z-z_r| = s'$ lies in the open set $ U' \subseteq U$ for all $s'$ close to $s$, 
and $v = \log |f/R| $ is harmonic on $U'$. Thus $J(s)$ exists and is given by 
\begin{eqnarray*}
J(s) &=& \frac1{2 \pi} \int_0^{2 \pi} \, \frac{\partial v }{ \partial \tau } \, ( z_r + s e^{i\theta} ) \, d \theta 
= 
\frac1{2 \pi} \int_0^{2 \pi} \, \frac{\partial \log |f |}{ \partial \tau } \, ( z_r + s e^{i\theta} )
 \, d \theta\\
&=& 
\frac1{2 \pi} \int_0^{2 \pi} \, \frac{\partial \arg f}{ \partial \theta  } \, ( z_r + s e^{i\theta} )
   \, d \theta = 
n(s, 1/F) - n(s, F) , 
\end{eqnarray*}
using the standard notation of Nevanlinna theory, where $ F(z) = f(z_r + z)$. 
Now any component $K'$ of $\D \setminus U'$ which meets 
$D(z_r, 4 \sigma) $ is contained in $D(z_r, s)$, and in particular this is true for $K' = K$.  
The boundary $\Gamma'$ of each such component $K'$ is a simple closed curve on which $|f| = R'$ and $\log f$ is locally univalent and
in particular sense preserving, 
by the choice of $R'$.
Since $\log (f/R')$ maps points lying just outside $K'$ into the right half plane it follows that,
as $z$ describes $\Gamma'$ once counter-clockwise, 
$\arg f(z)$ must increase. Thus the number of zeros of $f$ in each such $K'$ is at least one more than the number of poles. Since $f$ has neither
zeros nor poles in $U'$, it follows that
$J(s) \geq 1$ for all $s \in T$, which gives (\ref{atleast1}). 

Now (\ref{16}) and (\ref{atleast1}) deliver,  
as $s \to 2048 \sigma$ from below, 
$$
\log \frac{s}{1024 \sigma } \leq \int_{1024\sigma}^s \frac{J( t) \, dt}{t} = I(s) - I( 1024 \sigma ) \leq 
\frac1{2 \pi} \int_0^{2 \pi} u(z_r + s e^{i \theta} ) \, d \theta - u(z_r) \leq \phi(r). 
$$
Since $\phi(r)$ tends to $0$ as $r \to 1-$, this gives a contradiction if $r$ is close enough to $1$, and proves Lemma~\ref{disclem}. 
\hfill$\Box$
\vspace{.1in}

\begin{lem}
 \label{disclem2}
Let $r \in [r_0, 1) \setminus E$ be close to $1$. Then $f$ satisfies, for $z \in D(z_r, 2 \sigma)$,
 \begin{equation}
 \label{19}
\log \frac{f(z)}{f(z_r)} = a(r) \log \frac{z}{z_r} + g(z), \quad | g(z)| \leq 2 \phi (r) = o(1). 
\end{equation}
\end{lem}
\textit{Proof.} Here the logarithms are chosen so as to vanish at $z_r$. By Lemma \ref{disclem}, the function $v(z) $ 
is harmonic on $D(z_r, 4 \sigma)$ and equals  $ \log |f(z)/R|$ there, and on the same disc
$g(z)$ 
is analytic, with 
$${\rm Re} \, g(z) = u(z) \leq \phi(r) = o(1)$$
by (\ref{16}),  as well as $g(z_r) = 0$. Now for  $z \in D(z_r, 2 \sigma)$, applying
the Borel-Carath\'eodory inequality
leads to (\ref{19}), from which (\ref{con1}) follows at once. 
\hfill$\Box$
\vspace{.1in}

\begin{lem}
 \label{disclem3}
Fix $T \in (0, 2)$ and let $r \in [r_0, 1) \setminus E$ be close to $1$. Then $f$ satisfies (\ref{20}) for $z $ in $ D(z_r,  T \sigma)$.
In particular this holds for $T=1$. 
\end{lem}
\textit{Proof.} As in \cite{BRS}, this follows from (\ref{15}), (\ref{19}) and Cauchy's estimate for derivatives, which give
\begin{equation}
|g'(z)| \leq C \frac{ \phi(r)}{\sigma} \leq C \frac{ \phi(r)}{\varepsilon (r)} \leq C a(r)^{1-\beta}  
\label{gest}
\end{equation}
for $z \in D(z_r,  T \sigma)$, where $C$ is independent of $r$. 
The proof of Lemma \ref{disclem3} is complete, and so is that of Theorem \ref{thm1}. 
\hfill$\Box$
\vspace{.1in}

\section{The case of positive order: proof of Theorem \ref{thm2}} 

Assume the hypotheses of Theorem \ref{thm2}. 
Since Theorem \ref{thm1} has been proved, it 
suffices to show that the relation (\ref{20}) may be differentiated further to give an estimate for 
$f^{(q)}/f$ with $2 \leq q \leq M$. With  $(r_n)$  as in (\ref{a1b})  let $n$ be large and let
$r_n'$ and $r$ satisfy 
(\ref{a1c}). Then 
\begin{equation}
 \label{a3}
\frac1{a(r)^{1- \beta} } = \frac1{a(r)^{1- 2 \beta} a(r)^\beta }
\leq \frac1{a(r_n)^{1-2 \beta}a(r)^{ \beta} } = o \left( \frac{(1-r_n)^{1+\rho_0}}{a(r)^{ \beta} } \right) 
= o \left( \frac{1-r}{a(r)^{ \beta} } \right).
\end{equation}
It follows at once from (\ref{3a}) and (\ref{a3}) that 
\begin{equation*}
 \label{a(r)bigger}
2048 \sigma = \varepsilon (r) =  \frac{ 1}{  a(r)^{1-\beta}(\log a(r))^{1 +  \delta } } \quad \hbox{and} 
\quad \frac1{\varepsilon(r)} = o( a(r)) . 
\end{equation*}
In particular,  $\varepsilon (r) $ is as asserted in (\ref{a1d}), and Theorem \ref{thm2} follows from the next lemma. 

\begin{lem}
Suppose that $G \subseteq [r_0, 1) \setminus E$ is such that $\lim_{r \to 1-, r \in G} a(r) \varepsilon (r) = \infty $,
and let $M \in \N$. 
Then as $r \to 1-$ with $r \in G$ the function $f$ satisfies  (\ref{b1}) for $1 \leq q \leq M$ and 
$|z-z_r | < \varepsilon(r)/2048 $. 
\end{lem}
\textit{Proof.} 
It will  be proved by induction that (\ref{b1}) holds 
for $1 \leq q \leq M$ and $z \in D(z_r, (2- q/M) \sigma)$. Here the estimate (\ref{b1}) for $q=1$ follows from
Lemma \ref{disclem3} with $T = 2-1/M$. 
Assume next that $1 \leq q < M$ and that the assertion has been proved for $q$. Then (\ref{b1})  and Cauchy's estimate give
a positive constant $C$, independent of $r$ as long as $r \in G$, with 
$$
|L_q'(z)| \leq C \frac{a(r)^q }{\varepsilon(r)}  = o( a(r))^{q+1} 
$$
for $z \in D(z_r, (2- (q+1)/M) \sigma)$. Combining this estimate with the formula 
\begin{equation}
L_{q+1} = L_qL_1 + L_q' 
\label{Lq}
\end{equation}
completes the induction. 
\hfill$\Box$
\vspace{.1in}

\section{The case of zero order}

Assume throughout this section that $v$  has order $\rho (v) = 0$ in (\ref{a1}), let $E$ be as in Theorem
\ref{thm1}, and  denote positive constants by $C$. 

\begin{lem}
 \label{order0}
The functions $a(r)$ and $\varepsilon (r)$ satisfy 
$$a(r) \leq \left( \frac1{1-r} \right)^{1+o(1)} 
\quad \hbox{and} \quad a(r) \varepsilon (r) \to 0 \quad \hbox{as} \quad r \to 1-. $$
\end{lem}
\textit{Proof.} The first part is standard but the following details are included for completeness. As $r \to 1-$ set 
$s = 1 - (1-r)^2 $ and write  
$$
a(r) \log \frac{s}{r} \leq  \int_r^s \, \frac{a(t) \, dt}{t} \leq B(s) - B(r) \leq B(s) \leq (1-s)^{o(1)} \leq (1-r)^{o(1)} .
$$
Since
$$
\log \frac{s}{r} = \log \frac{2r-r^2}{r} = \log (1+1-r) \geq C (1-r) 
$$
as $r \to 1-$ the first assertion of the lemma follows, and so does the second, since (\ref{3a}) gives
$$
a(r) \varepsilon (r) \leq \frac{(1-r)a(r)^{1-\beta} }{2  (\log a(r))^{1 +  \delta } } \to 0. 
$$ 
\hfill$\Box$
\vspace{.1in}

In this zero order case further differentiation of (\ref{20}) need not lead to (\ref{b1}) for $q \geq 2$. Taking 
$1 < T < 2$ in Lemma~\ref{disclem3} and applying Cauchy's estimate to (\ref{20})  yields 
$$
\left| \frac{d}{dz} \left( \frac{f'(z)}{f(z)} \right) \right| \leq C \frac{a(r)}{\varepsilon(r)} ,
$$
for $z \in D(z_r,   \sigma)$,
but Lemma \ref{order0} shows that
the  upper bound arising here is not $o(a(r)^2)$. The example (\ref{exa1}) mentioned in the introduction shows that
this phenomenon is to be expected. 
However, an upper bound for $|f^{(q)}(z)/f(z)|$ is obtained as follows. 

\begin{thm}
 \label{order0a}
With the notation of Theorem \ref{thm1},  
assume that $v$  has order $\rho (v) = 0$ in (\ref{a1}), fix a positive integer $M$, and let  $r \in [r_0, 1) \setminus E$ be close to $1$. 
Then $f$ satisfies
\begin{equation}
 \label{order0est}
\left| \frac{f^{(q)}(z)}{f(z)} \right| \leq C a(r)  \left( \frac1{1-r} \right)^{(q-1)(1+\beta + o(1) )} 
\end{equation}
for $1 \leq q \leq M$ and $z \in D(z_r, (2- q/M) \sigma)$, where $\sigma = \varepsilon (r)/2048$.
\end{thm}
Since $\beta$ may be chosen arbitrarily small in (\ref{3}), the upper bound arising from Lemma \ref{order0} and Theorem \ref{order0a}
seems slightly stronger than that for analytic
functions of order zero in \cite{chyzh1,chyzh2}, but of course (\ref{order0est}) only holds near to the maximum
modulus.\\\\
\textit{Proof of Theorem \ref{order0a}.}
Define $L_q$ as in (\ref{b1}). 
For $q=1$ the asserted upper bound follows from Lemma \ref{disclem3} with $T = 2-1/M$. 
Assume next that $1 \leq q < M$ and that the asserted bound has been established for $q$.
Since $\beta $ is positive, Lemma \ref{order0} and (\ref{3a}) yield
$$
\varepsilon (r) = \frac{1-r}{2 a(r)^\beta (\log a(r))^{1 +  \delta } } \geq (1-r)^{1+\beta+o(1)} . 
$$
In conjunction with  Lemmas \ref{disclem3} and \ref{order0} and  Cauchy's estimate for derivatives this leads to
\begin{eqnarray*}
|L_q(z)L_1(z)| &\leq& C  a(r)^2  \left( \frac1{1-r} \right)^{(q-1)(1+\beta + o(1) )}  \\
&\leq& a(r) \left( \frac1{1-r} \right)^{(q-1)(1+\beta + o(1) ) + 1 + o(1) } \\
&\leq& a(r) \left( \frac1{1-r} \right)^{q(1+\beta + o(1) ) } 
\end{eqnarray*}
and 
\begin{eqnarray*}
|L_q'(z)| &\leq& C a(r) \,  \left( \frac1{1-r} \right)^{(q-1)(1+\beta + o(1) )} 
\,  \frac{1 }{\varepsilon(r)} \\
&\leq& C a(r) \,  \left( \frac1{1-r} \right)^{(q-1)(1+\beta + o(1) )}  \,  \left( \frac1{1-r} \right)^{1+\beta + o(1) }\\
&\leq &  a(r) \,  \left( \frac1{1-r} \right)^{q(1+\beta + o(1) ) } 
\end{eqnarray*}
for $z \in D(z_r, (2- (q+1)/M) \sigma)$. Using (\ref{Lq}), the theorem is proved by induction on $q$. 
\hfill$\Box$
\vspace{.1in}

\noindent
\textit{Acknowledgement.} The authors thank Walter Bergweiler, Peter Fenton and Janne Heittokangas for 
valuable correspondence.

\noindent
School of Mathematical Sciences, University of Nottingham, NG7 2RD.\\
jkl@maths.nott.ac.uk
\\\\
Department of Mathematics, 
Virginia Polytechnic Institute and State University, 
Blacksburg, VA 24061-0123. \\
rossij@vt.edu

{\footnotesize
\begin{thebibliography}{99}
\bibitem{Ber1}
W. Bergweiler, Proof of a conjecture
of Gross concerning fixpoints, \textit{Math. Zeit.} 204 (1990), 381-390.
\bibitem{Ber3}
W. Bergweiler, Periodic points of entire functions, proof of a conjecture
of Baker, \textit{Complex Variables Theory Appl.} 17 (1991), 57-72.
\bibitem{bergsizediscs}
W. Bergweiler, The size of Wiman-Valiron disks, 
\textit{Complex Var. Elliptic Equ.} 56 (2011), 13-33.
\bibitem{BRS}
W. Bergweiler, P.J. Rippon and G.M. Stallard,
Dynamics of meromorphic functions with direct or logarithmic
singularities, \textit{Proc. London Math. Soc.} 97 (2008), 368-400.
\bibitem{chyzh1}
I. Chyzhykov, G.G. Gundersen and J. Heittokangas, 
Linear differential equations and logarithmic derivative estimates,
\textit{Proc. London Math. Soc.} (3) 86 (2003), 735-754. 
\bibitem{chyzh2}
I. Chyzhykov, J. Heittokangas and J. R\"atty\"a, 
Sharp logarithmic derivative estimates with applications to ordinary differential equations in the unit disc,
\textit{J. Aust. Math. Soc.} 88 (2010), 145-167. 
\bibitem{Er0}
A. Eremenko,
On the iteration of entire functions,
\textit{Dynamical Systems and Ergodic Theory}, Banach Centre Publ. 23,
Polish Scientific Publishers, Warsaw 1989, 339-345.
\bibitem{FR}
P.C. Fenton and J. Rossi, ODEs and Wiman-Valiron theory in the unit disc, 
\textit{J. Math. Anal. Appl.} 367 (2010), 137-145.
\bibitem{Hayman}
W. Hayman, The local growth of power series: a survey of the Wiman-Valiron method,
\textit{Canad. Math. Bull.} 17 (1974) 317-358.
\bibitem{Lai1}
I. Laine, \textit{Nevanlinna theory and complex differential equations},
de Gruyter Studies in Math. 15,  Walter de Gruyter, Berlin/New York 1993.
\bibitem{Mac}
A.J. Macintyre, Wiman's method and the 'flat regions' of integral functions,
\textit{Quart. J. Math,. Oxford Ser.} 9 (1938) 81-88.
\bibitem{Tsuji}
M. Tsuji, \textit{Potential theory in modern function theory},
Maruzen, Tokyo, 1959.
\end{thebibliography}
}
\end{document}